\documentclass[a4paper]{article}

\pdfoutput=1

\usepackage{graphicx}
\usepackage{subcaption}
\usepackage{units}
\usepackage{amsmath}
\usepackage[subnum]{cases}
\usepackage{url}
\begin{document}

\title{Simultaneous optimisation of temperature and energy in linear energy system models}

\author{Patrik Sch\"onfeldt,
        Adrian Grimm,
        Bhawana Neupane,\\
        Herena Torio,
        Pedro Duran,
        Peter Klement,
        Benedikt Hanke,\\
        Karsten von Maydell, and
        Carsten Agert\\\\
        DLR Institute of Networked Energy Systems,\\
        Carl-von-Ossietzky-Str. 15, 26129 Oldenburg, Germany
}

\date{\today}

\maketitle

\begin{abstract}
Linear programming is used as a standard tool for optimising unit commitment or power flows in energy supply systems.
For heat supply systems, however, it faces a relevant limitation:
For them, energy yield depends on the output temperature,
thus both quantities would have to be optimised simultaneously
and the resulting problem is quadratic.
As a solution, we describe a method working with discrete temperature levels.
This paper presents mathematical models of various technologies
and displays their potential in a case study
focused on integrated residential heat and electricity supply.
It is shown that the technique yields reasonable results including the choice
of operational temperatures.
\end{abstract}

\section{Introduction}
\label{sec:intro}

Space heating accounts for approximately one third of the global
final energy consumption in both, the residential and the commercial building
sub-sectors.
Including hot water, the demand for low temperature heat was \unit[17.4]{PWh}
in 2010, which makes 53\% of the total final energy demand of
the worldwide building sector~\cite{IPCC_2014_Buildings}.
Still in 2019,
the share of renewable energy supply only meets 11\% of the global heat demand,
leading to a domination of fossil fuels in this sector,
contributing 40\% (\unit[13.3]{Gt}) of global CO2 emissions~\cite{IEA_renewables_2020}. 

This fact imposes a need for energy system optimisation in the heat sector.
As the efficiency of renewable heat sources such as solar thermal plants
or heat pumps shows a significant dependency on the temperature,
it should be considered in the optimisation.
However, while numerical descriptions including both,
energy and temperature, are common for
pure physical models (i.e.~\cite{FRANKE1997171,DECESAROOLIVESKI2003121,S0038092X01000366}),
cross-sectorial optimisation models typically lack according formulations
(cf.~\cite{PRINA2020109917,en12224298,FAZLOLLAHI2014648}).
A reason is that the commonly used technique of linear programming
does not allow straightforward implementations
that would require the product of optimisation variables.
To tackle this problem, the nonlinear problem can be modelled
(e.g.\cite{Mertz2017236,HUANG2019464,HOHMANN2019100177}),
but there are further options:

A possible solution is presented in~\cite{SCHUTZ201523}.
These authors subdivide the problem into two parts:
The storage is modeled by sub-volumes of constant size,
making temperature the optimisation variable.
Outside the storage, mass flows are optimised.
Alternatively, possible mass flows can be pre-calculated.
When deviations from this value are assumed to be
small, temperature and heat calculations can be decoupled~\cite{BAVIERE201869,giraud2017optimal}.
When iterating over such a linear model using different possible
mass flows, it will eventually converge.
An alternative approach works with discrete,
fixed values for either the mass flows~\cite{YOKOYAMA2017888},
or the temperature levels~\cite{FAZLOLLAHI2014648}.

The present work follows the latter approach
but introduces modifications to reduce the computational effort.
It also extends prior work using discrete temperature levels by
allowing for variable temperature levels at all parts of the model.
Also, an example study implementing the method using
\texttt{oemof.solph}~\cite{KRIEN2020100028} is presented.
In the study, the operation of a solar based residential energy supply
system is optimised either for economical costs or exergy.

\section{Mathematical models}
\label{sec:mathematical models}

The thermal energy of a volume \(V\) of matter at temperature
\(T\) with \([T] = \unit[1]{K}\) is
\begin{equation}
    E_\mathrm{th}(T, V) = \rho c_\mathrm{p} T V,
    \label{eq: E_thermal}
\end{equation}
where \(\rho\) is the density and \(c_\mathrm{p}\) is the thermal capacity
of the medium.
When there is a lower temperature level \(T_\mathrm{low}\) and
a higher one \(T_\mathrm{high}\),
it is convenient to also define a heat energy resulting from the temperature difference
\begin{equation}
    Q(T_\mathrm{s}) := E_\mathrm{th}(T_\mathrm{high} - T_\mathrm{low}).
    \label{eq: heat definition}
\end{equation}

Assuming a pipe diameter \(d\) the transferred heat then is
\begin{subequations}
\begin{align}
    \dot{Q}
        &= \rho c_\mathrm{p} \dot{V} \times (T_\mathrm{high} - T_\mathrm{low})\\
        &= \rho c_\mathrm{p} d v \times (T_\mathrm{high} - T_\mathrm{low}),
    \label{eq: heat transfer vT product}
\end{align}
\end{subequations}
where \(v\) marks the velocity of the matter at temperature
\(T_\mathrm{high}\).
Having in mind temperature-dependent losses and efficiencies,
it makes sense to optimise both, temperature and transferred volume
(\(\dot{V}\)) simultaneously.
However, as they are both factors in the same product,
linear programming does not allow to have variables for both
at the same time.
Thus it is not straightforward to optimise both using this optimisation technique.
As a workaround, we introduce \(k\) discrete temperature levels with
\begin{equation}
    T_\mathrm{high} \ge T_{n+1} > T_n \ge T_\mathrm{low}, \quad n \in 0 \dots k-1.
    \label{eq: temperature levels}
\end{equation}
For these fixed temperature levels, instead of Eq.~\eqref{eq: heat transfer vT product},
it is now possible to optimise
\begin{equation}
    \dot{Q} = \sum\limits_{n=0}^{k-1} \rho c_\mathrm{p} d \times v_n
              \times (T_n - T_\mathrm{low}),
    \label{eq: heat transfer vT sum}
\end{equation}
where special conditions might apply for the \(v_n\).
For example, if only one temperature level \(T_m(t)\) can be active at a time,
it would be \(v_n(t) = 0\; \forall n \neq m(t)\).

In the following, \(\dot{Q}_{\mathrm{in},n}\) identifies a heat flow to
the level at \(T_n\),
while \(\dot{Q}_{\mathrm{out},n}\) marks a heat flow from the level \(T_n\)
to another level at \(T_m \text{ with } n\neq m\).

\subsection{Heat sources}
\label{sec:heat sources mathematical models}

We assume that the temperature can always be reduced.
Thus, any heat source \(s\) supplying heat at \(T_n\) can supply
at least the same amount of heat at \(T_{n-1}\),
hence
\begin{equation}
    \dot{Q}_{s,n-1,\mathrm{max}} \ge \dot{Q}_{s,n,\mathrm{max}}.
    \label{eq: max heat supply temeprature}
\end{equation}
Also, heat is always gradually increased, from one temperature level to the next.
The heat entering the heat supply system \(\dot{Q}_\mathrm{in}\) at \(T_n\)
is composed of the supply coming from the external source at that level \(\dot{Q}_s(T_{n})\)
and heat taken from the lower level
\(\dot{Q}_\mathrm{out}(T_{n-1})\)
\begin{subequations}
\begin{equation}
\label{eq: temperature rise}
\begin{split}
    \dot{Q}_{\mathrm{in},n}
        &= (1-r_{n,n-1}) \times \dot{Q}_{\mathrm{out},n-1}\\
        &+ r_{n,n-1} \times \dot{Q}_{s,n},
\end{split}
\end{equation}
where
\begin{equation}
    r_{n,n-1} = \frac{T_{n}-T_{n-1}}{T_{n}-T_\mathrm{low}}.
\end{equation}
\end{subequations}
This guarantees that heat can not be transferred from a lower to a higher level without applying external work.

\subsubsection{Heat sources with constant efficiencies}

For heat sources with an efficiency that does not depend on the
target temperature, only the highest relevant temperature
output \(\dot{Q}_{s}(T_{s}) = \dot{Q}_{s,k-1}\) has to be modelled.
This is the case, e.g for boilers (gas, pallets, oil) or 
electric heating rods.
Equation~\ref{eq: max heat supply temeprature} is then implemented
using a possible lossless heat flow from
\(Q_n\) to \(Q_{n-1}\).
In combination with Eq.~\eqref{eq: temperature rise},
this yields
\begin{equation}
    \dot{Q}_s = \sum_{n} \dot{Q}_{\mathrm{in},n},
\end{equation}
which is the expected result.

\subsubsection{Heat sources with temperature-dependent efficiencies}
\label{sec: temperature-dependent heat sources}

Especially heat pumps or solar thermal plants --
heat sources that play a large role in renewable heat supply --
have temperature-dependent efficiencies \(\eta_s(T)\).
Having predefined temperature levels allows to pre-calculate these for all \(T_n\).
This way,
a maximum supply is modeled for every possible combination
of temperature level and supply technology,
\begin{equation}
    \dot{Q}_{s,n} \le \dot{Q}_{s,n,\mathrm{max}}\quad \forall\; n > 0.
\end{equation}
To avoid increased production of heat due to the increased number of sources,
binary status variables \(s_{s,n} \in 0, 1\) can be added,
with\footnote{Constraint \texttt{max\_active\_flow\_count} in oemof.solph.}
\begin{subequations}
\label{eq: sources exclusive status}
\begin{align}
    \sum_{n} s_{s,n} &\le 1 \quad \forall\; s\\
     \dot{Q}_{s,n} &\le s_{s,n} \times \dot{Q}_{s,n,\mathrm{max}} \quad \forall\; s, n.
\end{align}
\end{subequations}
While this approach is realistic,
it has to be considered that binary variables may introduce performance issues.
So, as an alternative to Eq.~\eqref{eq: sources exclusive status},
partial use at all levels can be allowed by introducing 
non-exclusive status variables \(q_{s,n} \ge 0\) with
\begin{subequations}
\label{eq: sources combined status}
\begin{align}
    \sum_{n} q_{s,n} &\le 1 \quad \forall\; s\\
     q_{s,n} \times \dot{Q}_{s,n,\mathrm{max}} &= \dot{Q}_{s,n} \quad \forall\; s, n.
\end{align}
\end{subequations}
If storage is in place, this can be justified by the presence of finite time steps
in the optimisation:
Even if the temperature levels have to be served exclusively at a certain time,
the level might be switched within one time step.

\subsection{Heat storage}
\label{sec:heat storage mathematical models}

Analogue to~\cite{FAZLOLLAHI2014648},
the total storage volume is split up into \(k\) sub-volumes
to map the temperature levels.
These form a stepped temperature-distribution model~\cite[Sec. 4.2.]{S0038092X01000366}
\begin{equation}
    V = \sum_{n=0}^{k-1} V_n(t)
    \label{eq: total storage volume}
\end{equation}
each with a homogeneous temperature \(T_n\)
that is chosen to be constant in time.
To reduce complexity, we exclude temperature inversion.
This is already implied by Eq.~\eqref{eq: temperature levels}.
Further, we assume the tank to be cylindrical.
Thus the volume is
\begin{equation}
    V_n(t) = h_n(t) \times \pi r^2
\end{equation}
and the surfaces are
\begin{numcases}
{A_n =\label{eq: moving boundaries - surfaces}}
2 \pi r h_n             & \, \(0 < n < k-1\) \\
2 \pi r h_n + \pi r^2    & \, \(n = 0 \lor n = k-1\),\label{eq: moving boundaries - surfaces top}
\end{numcases}

Although \(\sum \dot{V}_n = 0\),
by setting the lowest temperature level \(T_0 = T_\mathrm{i}\),
losses do not mean a heat exchange in the sense of \(\dot{Q}_n > 0\),
as for losses
\begin{subequations}
\begin{equation}
    \dot{V}_n \le 0 \quad \forall\; n > 1
\end{equation}
and due to the choice of \(T_0\)
\begin{equation}
    Q(V_0) \equiv 0.
\end{equation}
\end{subequations}

When setting \(T_0 = T_\mathrm{i}\) to the ambient temperature \(T_\mathrm{A}\),
the heat loss is
\begin{subequations}
\label{eq: moving boundary Qdot}
\begin{align}
    \dot{Q}_n
        &= (T_n - T_\mathrm{A}) \times \dot{V}_n \rho c_\mathrm{p}\\
        &= -\frac{\lambda}{d_\mathrm{iso}}(T_n-T_\mathrm{A})
            \times A_n(t).
\end{align}
\label{eq: moving boundary Qdot T}
\end{subequations}
Using Eq.~\eqref{eq: moving boundaries - surfaces},
\(h_\mathrm{n,0} = h_\mathrm{n}(t=0)\)
we get
\begin{numcases}
{h_\mathrm{n}(t) = \label{eq: two T levels h(t)}}
h_\mathrm{n,0} \times e^{-t/\tau_{\mathrm{MB}}}             & \, \(0 < n < k-1\) \\
h_\mathrm{n,0} \times e^{-t/\tau_{\mathrm{MB}}} + \pi \frac{r}{2}    & \, \(n = k-1\)
\end{numcases}
where the second term accounts for losses at the top side.
As the sub-volumes are independent,
they might be implemented as one individual storage for each
temperature level~\cite{MSc_Bhawana},
as long as Eq.~\eqref{eq: total storage volume} is preserved
by neglecting the \(n\) surfaces \(A_{0,n}\) of the lowest level, which is now present multiple times.

If the ambient temperature \(T_\mathrm{A}(t)\) is variable,
the stored heat according to Eq.~\eqref{eq: heat definition}
can change, even if nothing changes in the storage.
To account for this, virtual heat flows
\begin{equation}
    \dot{Q}_{\mathrm{A},n} = -\dot{T} V_n \rho c_\mathrm{p}
\end{equation}
have to be added for all sub-volumes \(V_n \text{ with } n > 0\)
that modify the stored energy accordingly.
Alternatively, a constant \(T_0 \approx T_\mathrm{A}(t)\) can be chosen,
if the losses for \(n=0\) remain small compared to those at
\(T_n\) with \(n > 0\).
This argument remains true if a representation based on multiple
individual storage tanks for the temperature levels is picked.
However, the influence of the 0th step on the losses will be overestimated
as there is an \(A_{n,0}\; \forall\; n > 0\).
In total
\begin{subequations}
\begin{align}
    \sum_{n=1}^{k} (A_{n,0})
    &= (k-1) \times A - \sum_{n=1}^{k} A_n\\
    &= (k-2) \times A + A_0,
\end{align}
\end{subequations}
meaning the surface of the 0th level is overestimated
by an offset proportional to the total surface of the storage.

Comparing our approach with the one of~\cite{FAZLOLLAHI2014648},
there are some significant differences:
While we model losses as virtual heat flows directly to the lowest level,
they connect adjacent sub-volumes~\cite{BECKER2012415}.
While the latter is more realistic,
our approach allows for dropping the binary variable
that indicates whether or not a temperature level is active.
This fact alone typically comes with an advantage in computational performance~\cite{1664974,ALEMANY2018429}.
Also, it is the key to further simplify the formulation by allowing alternative representation using \(k-1\) storage tanks with just two layers, each.

\subsection{Heat demands}
\label{sec:heat demands mathematical models}

Heat demands \(\dot{Q}_{d}\) are modeled by a reduction of the temperature,
regardless if that reduction really happens or if fresh, cold water is replacing hot water in the tank (i.e. when showering).
Analogue to Eq.~\eqref{eq: temperature rise}, we have
\begin{subequations}
\label{eq: temperature reduce}
\begin{equation}
\begin{split}
    \dot{Q}_{\mathrm{out},n}
        &= r_{n,m} \times \dot{Q}_{d,n}\\
         & + (1 - r_{n,m}) \times \dot{Q}_{\mathrm{in},m},
\end{split}
\end{equation}
where
\begin{equation}
    r_{n,m} = \frac{T_{n}-T_{m}}{T_{n}-T_\mathrm{low}}.
\end{equation}
\end{subequations}
So, just in the special case of \(T_m = T_\mathrm{low}\)
all heat taken at level \(T_n\) fulfills
the the demand at that level.
Typically it partly increases the heat at the level
\(T_m < T_n\).

\section{Case study}

\subsection{Energy system layout}
\label{sec:example energy system layout}

\begin{figure}[tbh]
    \centering
    \includegraphics[width=0.95\linewidth]{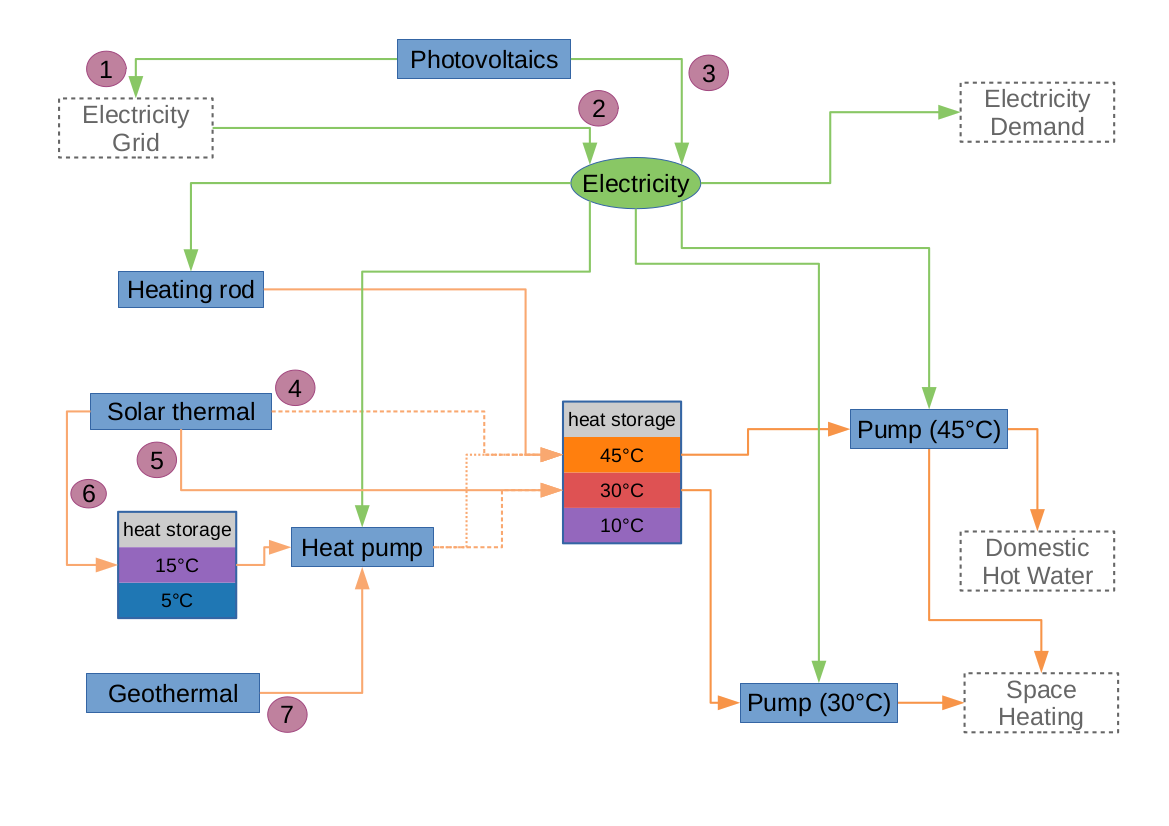}
    \caption{Energy system graph of example scenario.
    Lines signify possible energy flows between components,
    dashed or dotted lines have no special meaning,
    the individual styles just makes them distinguishable.
    The numbers mark energy flows that have an attached cost
    (see Tab.~\ref{table:costs_ren_data}).}
    \label{fig: scenario graph}
\end{figure}

As an example, we optimise the control strategy
of an integrated energy supply system of a residential block
using the presented techniques.
The parameters have been chosen to be realistic for
the quarter ``Helleheide''~\cite{helleheide} located
in north-western Germany,
which is being developed in the framework of the research project “Energetisches Nachbarschaftsquartier Fliegerhorst Oldenburg” (short ENaQ).
The used energy system model consists of a number of energy supply and conversion technologies
as well as two heat storage tanks.
A graph representing this energy system is depicted
in Fig.~\ref{fig: scenario graph}.

We assume about 140 residential units, resulting in a demand for electricity of \(E_\mathrm{el} = \unitfrac[376]{MWh}{a}\)
and heat demands of \(Q_\mathrm{SH} = \unitfrac[402]{MWh}{a}\)
at at least \(\unit[30]{^\circ C}\) for space heating
and \(Q_\mathrm{DHW} = \unitfrac[336]{MWh}{a}\)
at \(\unit[45]{^\circ C}\) for domestic hot water.
Besides economical costs (in EUR),
we optimise for an exergy rating \(X(T_\mathrm{A}(t))\)
that varies with the ambient temperature.
Electrical power to pump the heat to the demand is assumed to be
\(\unit[1]{\%}\) of the transported heat for the higher and
\(\unit[2]{\%}\) of the transported heat for the lower
temperature level.

On the supply side, there are a PV power plant
(\(P_\mathrm{peak} = \unit[150]{kW}\)),
solar thermal collectors (\(A = \unit[1050]{m^2}\)),
and a geothermal collector
(\(P_\mathrm{th,max} = \unit[100]{kW}, T = T_\mathrm{soil}(t)\)).
All three show weather dependent characteristics.
The solar thermal collector is implemented using
the technique described in
Sec.~\ref{sec: temperature-dependent heat sources},
so it reflects the trade-of between energy gain and
achievable temperature.
The same is true for the heat pump
(\(P_\mathrm{th} = \unit[280]{kW}\)),
which can be used to raise the temperature from either
\(T_\mathrm{soil}(t)\) or \(\unit[15]{^\circ C}\) to
either \(\unit[30]{^\circ C}\) or \(\unit[45]{^\circ C}\).
For this,
coefficients of performance (COPs) of \(\unit[30]{\%}\) of the Carnot COP \(1/\eta_\mathrm{c}\),
see Eq.~\eqref{eq: carnot efficiency}, are considered.

Further, there are a cold and a warm heat storage
(\(V = \unit[50]{m^3}\) each).
These are modeled according to
Sec.~\ref{sec:heat storage mathematical models}.
The cold storage solely serves as a source for the heat pump
and uses \(T_0 = \unit[5]{^\circ C}\).
The warm storage is used to fulfill the demands and has
\(T_0 = \unit[10]{^\circ C}\),
which is considered the temperature of fresh water.
Both storage tanks can get energy from the solar thermal collector.
Additionally, the warm storage is fed by the heat pump
or an electric heating rod \(\eta = \unit[95]{\%}\) serving as a backup.
Note that the energy system layout has been chosen to see the applicability and plausibility of the method on a realistic heat energy supply system. 
It does not necessarily represent an energy system of optimal performance though.

\subsection{Input data}

% cost data in renewable-based scenario
\begin{table}
    [tbh]
    \centering
    \begin{tabular}{ |l|r|r| } % 4 columns
         \hline % column headers
             no. & Prices & Exergy \\
                 & (\unitfrac{EUR}{MWh}) & (\unitfrac{MWh}{MWh}) \\
         \hline % column contents
            1 & -72.90  & 0 \\
            2 & 34.07 \dots 280.65 & 1.224 \dots 1.704 \\
            3 & 67.56 & 1.0 \\
            4 & 0 & 0.000 \dots 0.080 \\
            5 & 0 & 0.000 \dots 0.126 \\
            6 & 0 & 0.047 \dots 0.167 \\
            7 & 0 & 0.000 \dots 0.057 \\
         \hline
    \end{tabular}
    \caption{
        Costs in optimisation of the model (Fig.~\ref{fig: scenario graph}).
        The dots signify a range in which the value varies
        within the year 2017,
        which was taken as an example in the study.
    }
    \label{table:costs_ren_data}
\end{table}

Input data consists of hourly time series for prices, weather,
and grid electricity for the year 2017.
An overview is given in Tab.~\ref{table:costs_ren_data}.
Further, there are demand time series with the same resolution.
The time series for space heating has been prepared according to~\cite{en13112967} using \texttt{QuaSi}~\cite{QuaSi},
the ones for electricity and domestic hot water using the
\texttt{Load Profile Generator}~\cite{PFLUGRADT2017655}.

\begin{figure}[tbh]
    \centering
    \includegraphics[width=\linewidth]{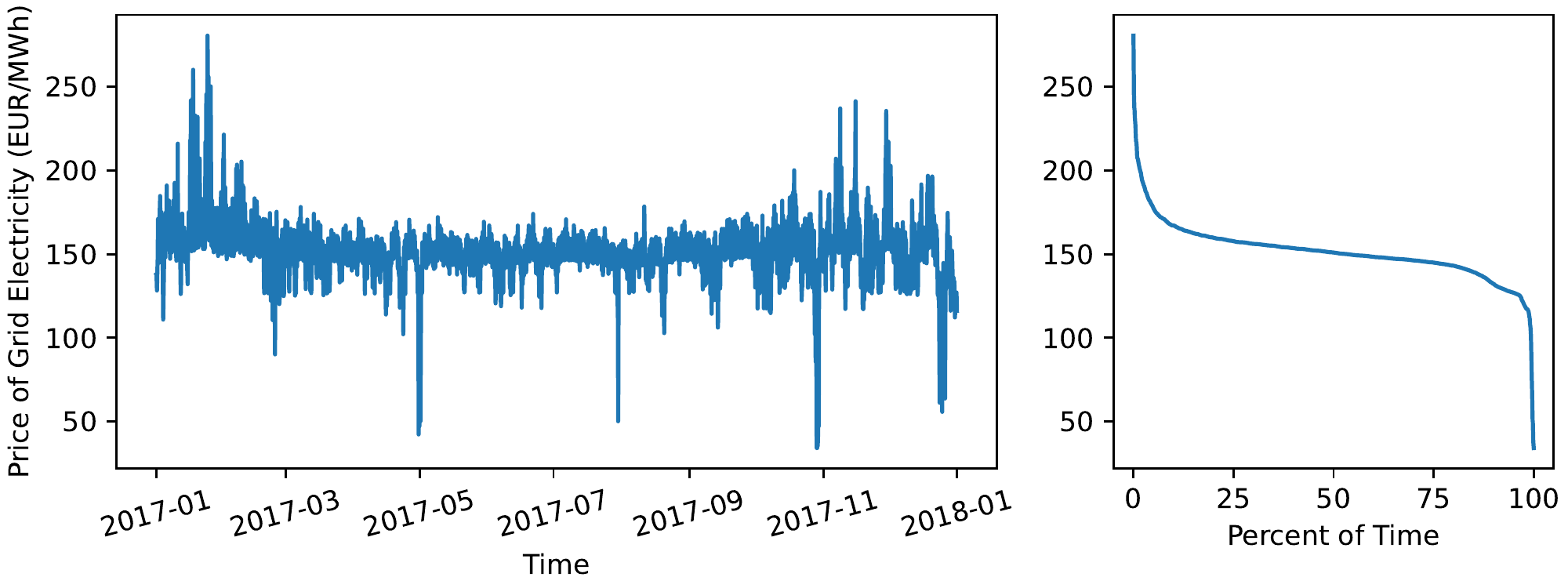}
    \caption{hour-specific costs for imports of electricity in terms of prices over time for two exemplary weeks of the year}
    \label{fig: costs_electricity_price}
\end{figure}

The applicable prices include (see~\cite[p.~41]{MSc_Adrian} for details)
\begin{enumerate}
    \item funding for pv feed in to the grid
    \item the day ahead price~\cite{EntsoeTransparancyPlattform} plus all applicable taxes and levies, and
    \item renewable energy levy according to German EEG.
\end{enumerate}
Figure~\ref{fig: costs_electricity_price} shows
the variations over the year and a price-duration curve.

% figure: costs for grid imports in terms of prices

The exergy of electricity from the grid is based on \cite{en13112967}.
There, the emissions from electricity in the grid are calculated
as an hour-specific average and dependent on the power plants operating
in a particular region within the European transmission grid.
In the same way, we calculate the share of fossil and renewable energy
feed-ins for each hour in 2017,
based on data of the entso-e~\cite{EntsoeTransparancyPlattform}
for the region \textit{DE-AT-LU}.
% figure: costs for grid imports in terms of exergy
\begin{figure}[tbh]
    \centering
    \includegraphics[width=\linewidth]{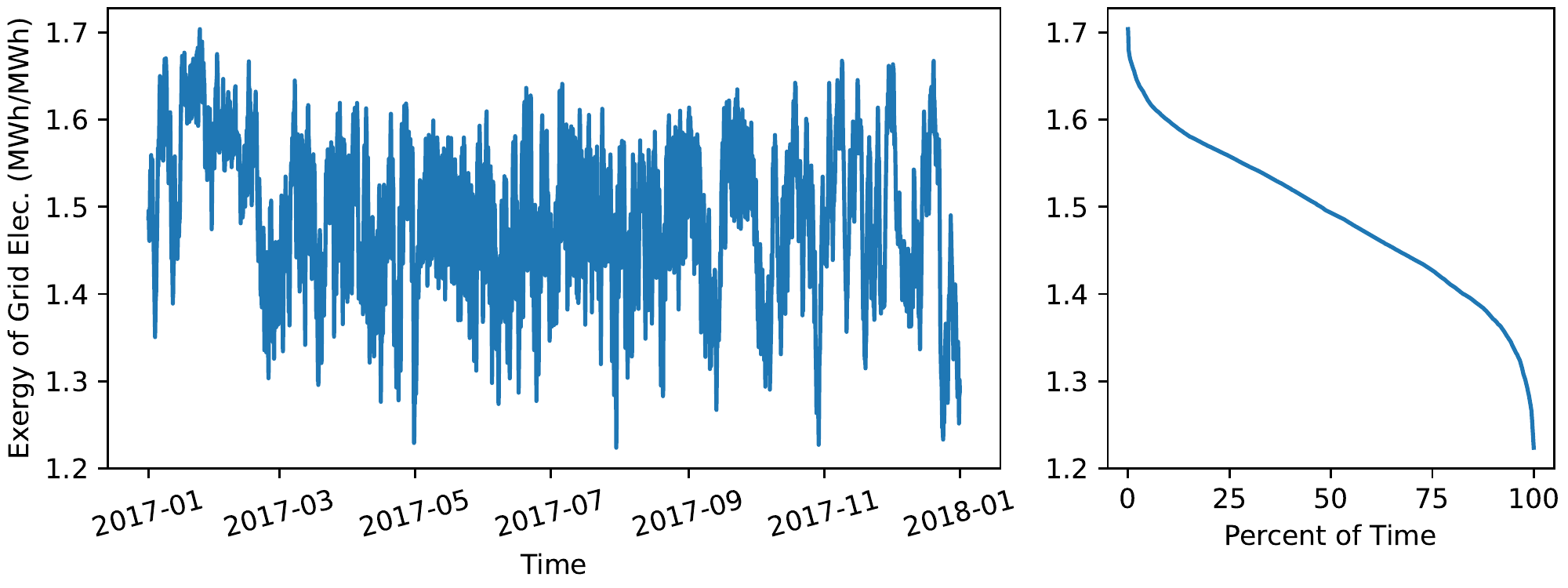}
    \caption{hour-specific costs for imports of electricity in terms of relative exergy-input over time for two exemplary weeks of the year}
    \label{fig: costs_electricity_exergy}
\end{figure}
For fossil-based energy, we assume \unitfrac[1.8]{MWh}{MWh},
which is the primary energy factors of the power plants
according to~\cite{SachstandPrimarenergiefaktoren2016},
as a best available guess for the relative exergy-input.
For renewable-based grid electricity, we set \unitfrac[1.05]{MWh}{MWh},
considering a 5\% loss in the grid,
local PV electricity has a value of \unitfrac[1.0]{MWh}{MWh}.
The combination of these data leads to the minimum and maximum value
for the exergy of the electricity from the grid shown in Table~\ref{table:costs_ren_data},
The values over the year and a duration curve are displayed in
Fig.~\ref{fig: costs_electricity_exergy}.
As the temperatures are assumed to be constant,
the other exergy ratings are calculated using the 
Carnot efficiency
\begin{equation}
    \eta_\mathrm{c} = 1 - \frac{T_\mathrm{low}}{T_\mathrm{high}},
    \label{eq: carnot efficiency}
\end{equation}
where \(T_\mathrm{low}\) and \(T_\mathrm{high}\)
are the minimum and the maximum
of the air temperature the temperature valid for the technology, respectively.
So, the weights are 
\begin{enumerate}\setcounter{enumi}{3}
    \item \(\eta_{15}(t) = 1 - \min(\unit[288.15]{K},T_\mathrm{A})
                        / \max(\unit[288.15]{K},T_\mathrm{A})\)
    \item \(\eta_{30}(t) = 1 - \min(\unit[303.15]{K},T_\mathrm{A})
                        / \max(\unit[303.15]{K},T_\mathrm{A})\)
    \item \(\eta_{45}(t) = 1 - \min(\unit[318.15]{K},T_\mathrm{A})
                        / \max(\unit[318.15]{K},T_\mathrm{A})\)
    \item \(\eta_\mathrm{soil}(t) = 1 - \min(T_\mathrm{soil},T_\mathrm{A})
                                      / \max(T_\mathrm{soil},T_\mathrm{A})\)
\end{enumerate}
for the corresponding energy flows.

\subsection{Optimisation results}

\begin{figure}[p]
    \centering
    \begin{subfigure}[b]{\textwidth}
        \centering
        \includegraphics[width=0.95\textwidth]{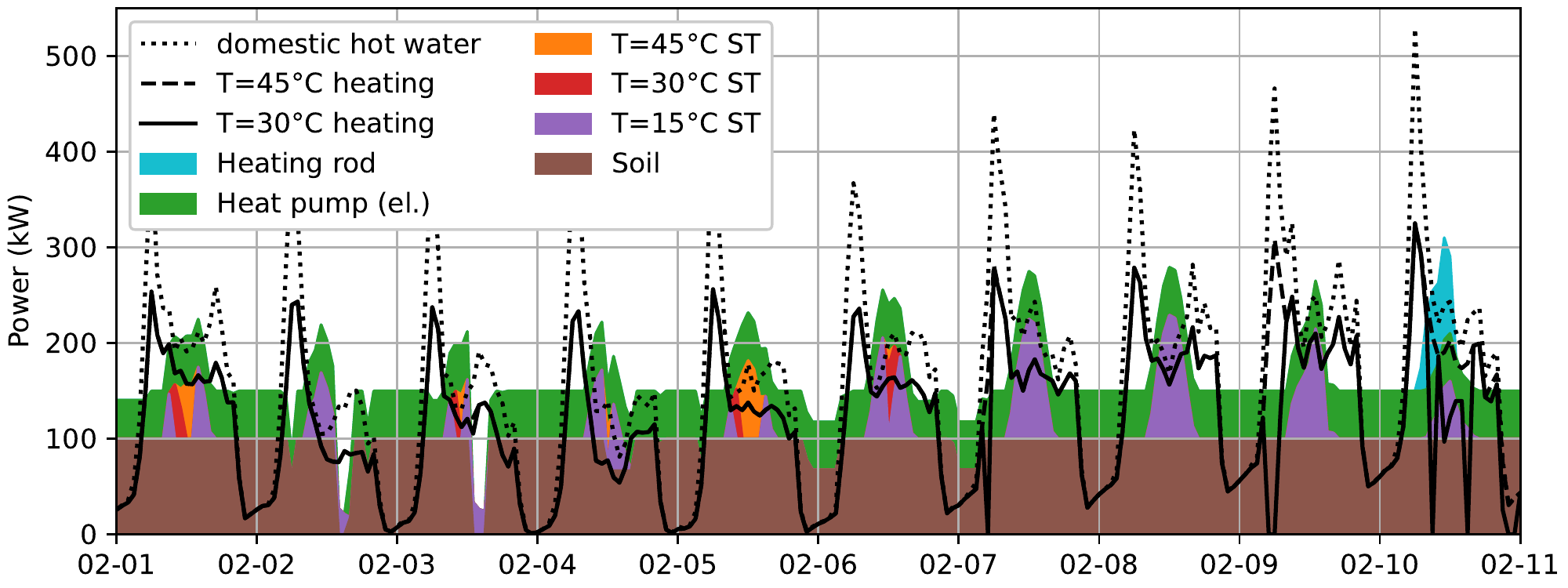}
        \caption{Hourly heat balance in the exergy-optimised case.}
        \label{fig: balance_optimised4exergy}
    \end{subfigure}
    \\
        \begin{subfigure}[b]{\textwidth}
        \centering
        \includegraphics[width=0.95\textwidth]{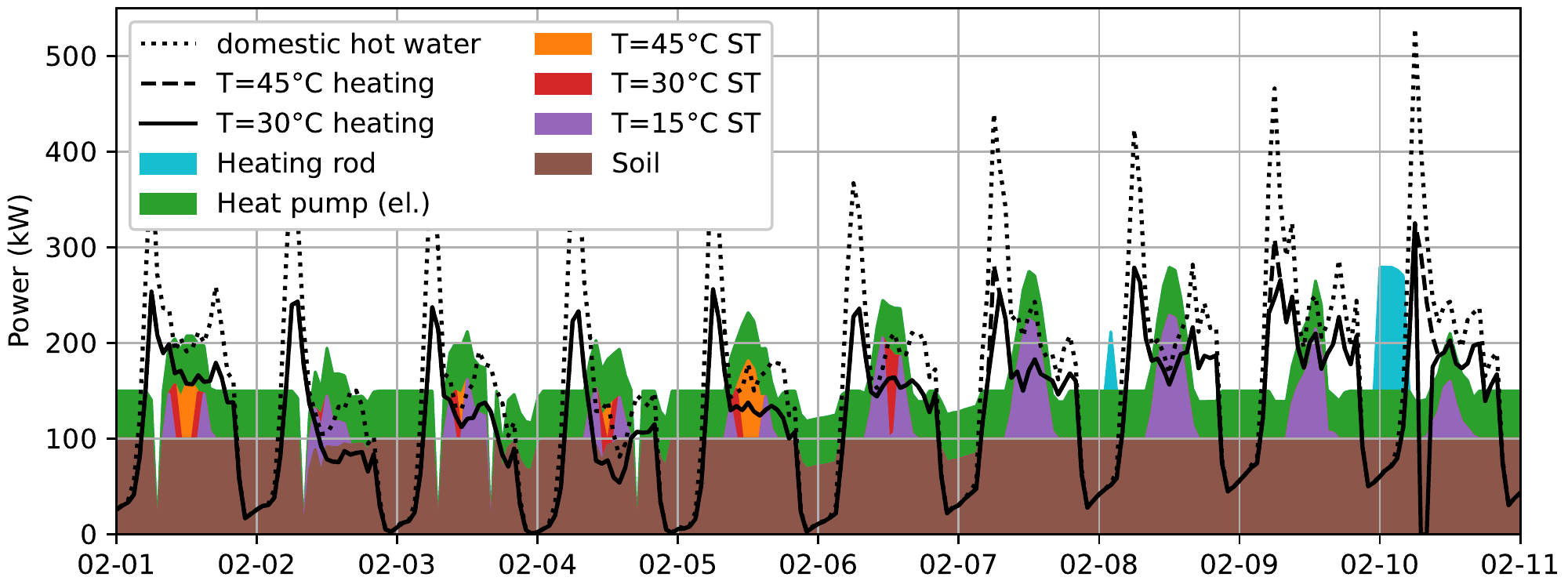}
        \caption{Hourly heat balance in the price-optimised case.}
        \label{fig: balance_optimised4prices}
    \end{subfigure}
    \\
    \begin{subfigure}[b]{0.48\textwidth}
        \includegraphics[width=\textwidth]{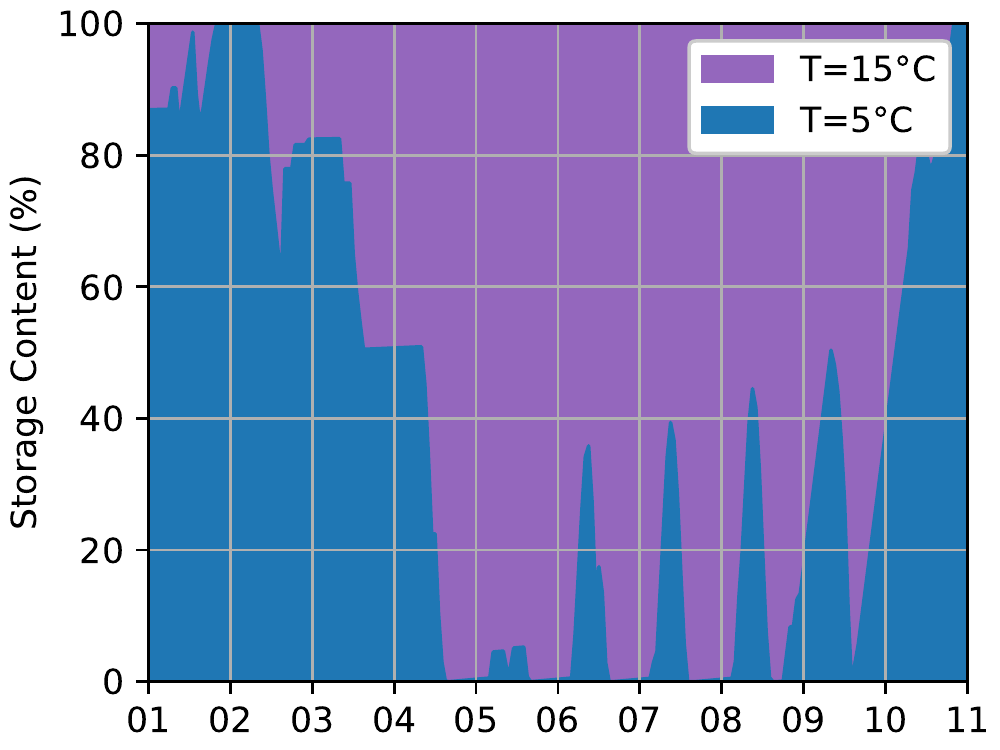}
        \caption{Cold storage, exergy-optimised}
        \label{fig: storage_1_exergy}
    \end{subfigure}
    \hfill
    \begin{subfigure}[b]{0.48\textwidth}
        \includegraphics[width=\textwidth]{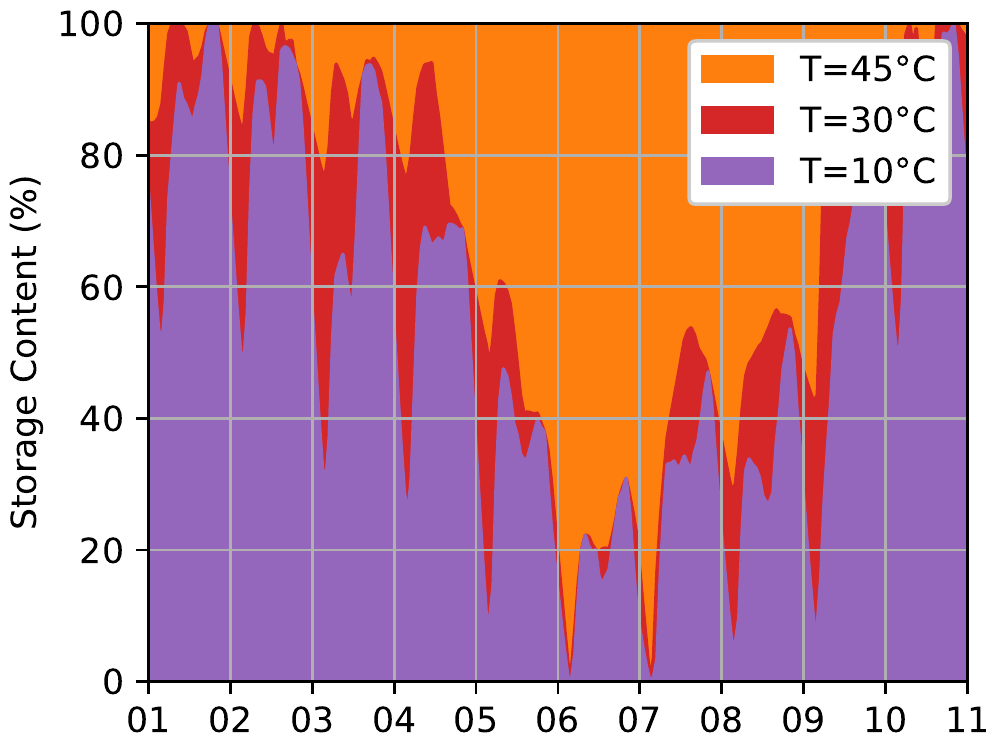}
        \caption{Warm storage, exergy-optimised}
        \label{fig: storage_2_exergy}
    \end{subfigure}
    \\
    \begin{subfigure}[b]{0.48\textwidth}
        \includegraphics[width=\textwidth]{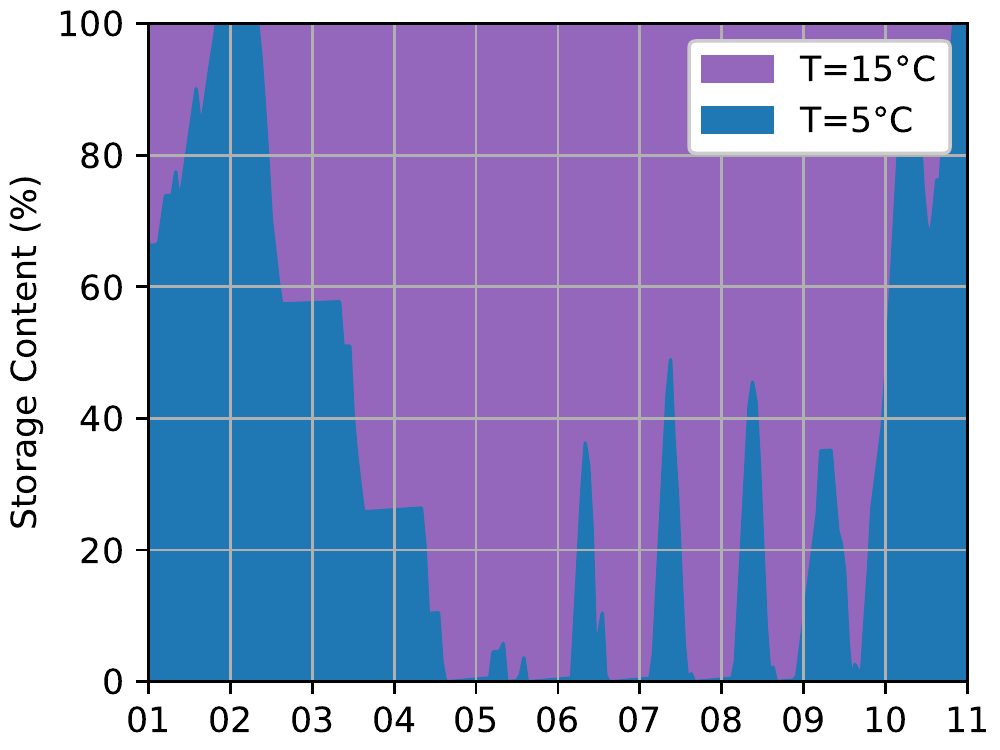}
        \caption{Cold storage, price-optimised}
        \label{fig: storage_1_prices}
    \end{subfigure}
    \hfill
    \begin{subfigure}[b]{0.48\textwidth}
        \includegraphics[width=\textwidth]{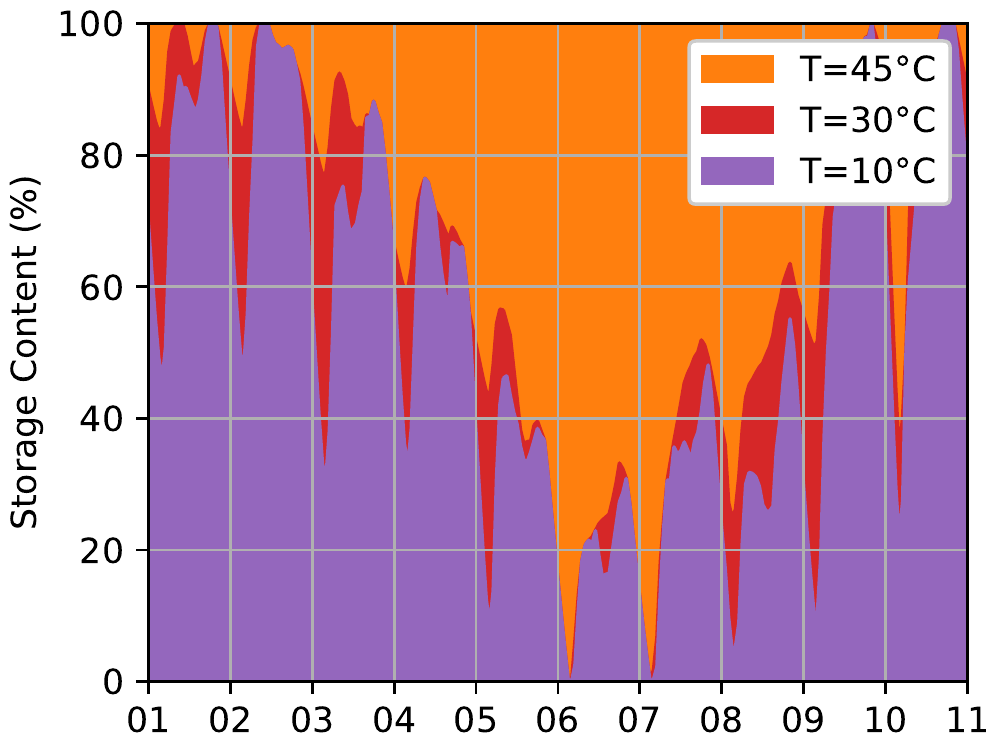}
        \caption{Warm storage, price-optimised}
        \label{fig: storage_2_prices}
    \end{subfigure}
    \caption{Hourly heat balance and storage contents
        for the first half of February.
        For \subref{fig: balance_optimised4exergy},
        \subref{fig: balance_optimised4prices},
        the lines mark demands, the colored areas denote supply.
        Panels \subref{fig: storage_1_exergy}
        through~\subref{fig: storage_2_prices}
        indicate usage of the volume.}
    \label{fig: optimisation results}
\end{figure}

The optimisation is run for the whole year 2017,
for both, prices and exergy.
For the first ten days of February,
the heat energy balance (heat being produced and used)
and the storage content (in relative volume) are displayed
in Fig.~\ref{fig: optimisation results}.
The values for the variables are displayed in
Table~\ref{tab: optimisation results}.

In Figs.~\ref{fig: balance_optimised4exergy}
and~\subref{fig: balance_optimised4prices},
it can be seen that in this period of time,
the soil is the main source of heat.
Additional energy is added using the heat pump when raising the temperature level.
The COPs depend on \(T_\mathrm{soil}\), they are
\(3.44 < 1/\eta_{30} < 7.00\) and \(2.30 < 1/\eta_{45} < 3.41\), respectively.
In the last days, the solar thermal panels provide only low-temperature
heat at the level of \(T=\unit[15]{^\circ C}\),
on the earlier days, higher temperatures are provided.
The utilisation of the storage tanks is similar for both scenarios
(see Fig~\ref{fig: storage_1_exergy}
through~\subref{fig: storage_2_prices}), as well:
In the beginning of the chosen period,
the storage tanks are almost completely at the lowest possible temperature level.
Until the sixth to seventh day, the storage tanks are then heated up.
The next few days, the sun can fully recover the cold storage by day,
but the warm storage is gradually emptied regardless of the
optimisation goal.
In the end of the period, both storage tanks are empty, again.
Also note that the level of \(\unit[45]{^\circ C}\) is occasionally used for heating.
This is due to the way that auxiliary energy for pumping is considered in this specific case study.
When heat is provided by the heating rod, i.e. at the 10th of February, the overall efficiency of this technology plus pumping favors the higher level.
Heating efficiency is considered to be constant for the two possible levels and electricity for pumping can be saved by providing less volume at the higher
temperature.

\begin{table}
    \centering
    \begin{tabular}{ |c|l|r|r| }
        \hline % column headers
        Quantity & Unit & Price-optimal & Exergy-optimal \\
        \hline % column contents
         Operational costs & \unitfrac{EUR}{MWh} & 53.68 & 54.16\\
         Exergy per Energy & \unitfrac{MWh}{MWh} & 0.617 & 0.611\\
         Specific emissions & \unitfrac{kg}{MWh} & 168 & 167\\
         Solar heat (\(\unit[15]{^\circ C}\)) & \unitfrac{MWh}{a} & 25.7  & 27.6\\
         Solar heat (\(\unit[30]{^\circ C}\)) & \unitfrac{MWh}{a} & 14.9 & 75.1\\
         Solar heat (\(\unit[45]{^\circ C}\)) & \unitfrac{MWh}{a} & 270 & 210\\
         Pumping electricity & \unitfrac{MWh}{a} & 10.6 & 11.1\\
        \hline
    \end{tabular}
    \caption{Resulting values for both optimisation cases.}
    \label{tab: optimisation results}
\end{table}

When optimising for prices or exergy,
the overall trend of the heat balance is very similar.
This is due to the fact that the chosen energy system design
does not leave much flexibility to alter the control strategy
when the other optimisation goal is followed.
However, there are significant differences:
In particular, in the price-optimised case,
the heating rod is operated at nights,
where the prices are lower.
In the exergy-optimised case, it is operated by day,
as the share of renewable energies tends to be higher by day,
resulting in a lower exergy-rating.
This difference can also be seen in the content of the warm storage,
which is used to store \(\unit[45]{^\circ C}\) at the last day
for the price-optimal but not for the exergy-optimal case.
Another difference is visible for the 2nd and 3rd day:
In the exergy-optimal case,
heat is preferably drawn from the cold storage,
while the geothermal source is preferred for the price-optimal case.
Considering the similarities in the temporal profiles,
it is no surprise that the final value for the price or the exergy,
only increases by 1\%, when the other variable is optimised 
(see Table.~\ref{tab: optimisation results}).
Over the full year, the most significant difference is
seen in the operation of the solar plant.
The price-optimal solution prefers the level at
\(\unit[45]{^\circ C}\) more often,
while the exergy-optimal one prefers \(\unit[30]{^\circ C}\).
The latter comes at the cost of a \unit[4.7]{\%} increase
of needed pumping electricity.

\section{Conclusion}

This paper presents a method to simultaneously optimise both,
temperature and energy within energy systems, using linear programming.
It uses discrete temperature levels and can work without binary or integer variables.
Using the technique, an example study was conducted,
showing the potential of the method.

The study optimises unit commitment, including the choice of
operational temperatures, for time-dependent exergy ratings and prices.
Even though the example energy system does not allow for much flexibility,
differences between the two optimisation scenarios are visible.
In particular, differences have been shown in the operation of the heat pump,
in switching supply temperatures and in the timings for using the heating rod.

Concluding, the presented method is capable of optimising
unit commitment and operation of heat supply including the
choice of operational temperatures, using the well-established
method of linear programming.

\section*{Author contributions}
Conceptualization, methodology, validation: P.S., H.T., P.D., B.H.,  K.v.M., and C.A.;
software, formal analysis, investigation, data curation: P.S., A.G., P.D., and B.N.;
writing--original draft preparation: P.S. and A.G.;
writing--review and editing: P.S., A.G., B.N., H.T., P.D., P.K. and B.H.;
visualization: P.S. and A.G.;
supervision: P.K. and B.H.;
project administration: P.K.;
funding acquisition: P.K., B.H., K.v.M., and C.A.;
All authors have read and agreed to the published version of the manuscript.

\section*{Acknowledgements}
This work has been funded by the
Federal Ministry for Economic Affairs and Energy (BMWi) of Germany and
the Federal Ministry of Education and Research (BMBF) of Germany
(grant number 03SBE111).

\section*{Conflict of interest}
The authors declare that they have no conflict of interest.

\nocite{IPCC_2014}

\bibliographystyle{unsrt}       % APS-like style for physics
\bibliography{preprint.bib}   % name your BibTeX data base

\end{document}